# A two-stage stochastic MINLP model to design and operate a multi-energy microgrid by addressing carbon emission regulatory policies uncertainty

Handan Akülker[a,b*], Burak Alakent[a], and Erdal Aydin[c]

[a]*Department of Chemical Engineering, Bogazici University, Istanbul, 34342, Türkiye*
[b]*Department of Chemical Engineering, Ondokuz Mayis University, Samsun, 55139, Türkiye*
[c]*Department of Chemical and Biological Engineering, Koc University, Istanbul, 34457, Türkiye*

*Corresponding author: handan.akulker@omu.edu.tr*

## Abstract

This study suggests a novel two-stage Mixed-Integer Nonlinear Programming model considering uncertainty related to implementation of carbon dioxide emission regulatory policies, which are carbon trading and emission taxing and can change over the years, for the purpose of optimal equipment selection from candidate equipment to design, size and operate a multi-energy microgrid. The uncertain sources are air temperature, wind speed, solar radiation, carbon dioxide trading price or tax, and natural gas price. Candidate equipment are wind turbines, PV arrays, a biomass-fired generator, biomass combined cycles, combined heat and power generators, conventional generators, an electricity storage unit, integrated gasification combined cycles, a heat pump, and a power-to-synthetic natural gas (P2G) system. Three case studies are investigated. In the first case, the model selects the optimal equipment for meeting the electricity and heat demands only. In the second case, the optimal equipment selections are determined to couple with the P2G system to meet the electricity, heat, and natural gas demands. In the third case, the model selects the optimal equipment to run with sustainable energy generators: wind turbines and solar panels. The optimal selections are compared between deterministic and stochastic forms of the optimization models.

# 1. Introduction

Multi-energy microgrids (MEMGs) may utilize various energy sources such as solar radiation, wind, waste biomass, oil, natural gas, and coal. Wind turbines, photovoltaic panels, combined heat and power (CHP) generators, biomass-fired generators, batteries, and conventional power generators may be alternative equipment to design a MEMG for energy producers (Chen et al., 2021; Mansour-Saatloo et al., 2020). However, during the optimal design and scheduling, energy producers should consider some issues such as carbon dioxide emissions, natural gas shortages, and volatile resource prices.

Carbon dioxide emission is an active problem for many industries, and primarily for the energy industry. The governments have started setting policies for energy producers to eliminate carbon dioxide emissions. There exist two different policies, namely carbon emission taxing and carbon trading (cap-and-trade system). Carbon emission taxing dictates taxes per generated amount of carbon dioxide emission. On the other hand, carbon trading promotes selling or purchasing carbon dioxide emission limits. Carbon trading price or emission tax may vary from one country to another around the world (Savaresi, 2016; TheWorldBank, 2024).

Natural gas is one of the major source of energy producers [5]. Power-to-Gas systems comprise power-to-methane (or synthetic natural gas) and power-to-hydrogen systems. They both consume excess electricity of energy grids for water electrolysis to produce pure hydrogen. On the other hand, power-to-methane systems extract carbon dioxide in high purity from flue gases or air to feed with pure hydrogen to synthesize methane, and by upgrading it, synthetic natural gas (SNG) is obtained. Thus, power-to-synthetic natural gas systems offer double benefits by reducing carbon dioxide emissions and producing SNG to handle the Paris Agreement sanctions and the natural gas supply problem (Glenk and Reichelstein, 2019; Mansouri et al., 2021; Wang et al., 2019). Therefore, it is also an active research problem to account for synergies and integration/symbiosis opportunities between chemical production plants and energy generation systems.

In the literature, two different deterministic MINLP-based decision-making models including carbon emission taxing and trading were proposed to meet the electricity demand (Akulker and Aydin, 2023). In that study, the effects of carbon taxing and trading on the optimal equipment selections were compared. In addition, the impacts of different carbon dioxide trading or taxing price, and various carbon dioxide limits were investigated. The carbon trading price was assumed to be the same as the emission tax for fair comparison in that study. Candidate equipment are wind turbines farms, solar panel arrays, a biomass-fired generator, biomass combined cycles, CHP units, conventional generators, an electricity storage unit, integrated gasification combined cycles, and a power-to-natural gas (P2G) system. Later, in the second step, a deterministic MINLP model was suggested for optimal equipment selections to couple with a P2G system under carbon dioxide taxing to meet the electricity and natural gas demands (Akulker and Aydin, 2024). The optimal selections based on different natural gas prices and carbon dioxide taxes were compared. On the other hand, these studies did not take into account the inherent uncertain mechanisms of resources and carbon adjustment policies.

Accordingly, in this work, a novel two-stage stochastic MINLP model is developed to meet the electricity, heat, and gas demand to investigate the differences from the generated deterministic models in our previous studies (Akulker and Aydin, 2024, 2023). Moreover, candidate equipment modeling, except for the heat pump, is the same. A heat pump and a CHP unit are additionally introduced compared to to meet the heat demand. The uncertainty sources considered in this study are carbon dioxide regulation policies (carbon trading or emission taxing), weather conditions (solar radiation, air temperature, and wind speed), carbon dioxide taxing or trading price, and natural gas price. Three case studies are investigated. In the first case, the optimal equipment selection is determined to meet the electricity and heat demand. In the second case, the purpose is to optimally design and schedule a P2G-integrated MEMG by selecting equipment to meet the electricity, heat, and natural gas demands. The hourly sale of generated SNG to the national gas grid is kept constant. For the first and second cases, the results of stochastic models are compared with the results of deterministic forms of the models. Finally, in the third case, the optimization model selects the equipment to work with sustainable energy generators: wind turbines and solar panels.

**2. Related Work**

The optimal design and scheduling for MEMGs are mostly studied by performing Mixed-Integer Programming (MIP) (Bartolucci et al., 2022; Hurwitz et al., 2020; Lekvan et al., 2021). MIP studies can be classified as deterministic or stochastic. Deterministic MIP-based models underestimate the effects of uncertainties depending on market price, weather conditions, and demands (Han and Lee, 2021; Wang et al., 2022; Zhang et al., 2020). Instead, stochastic MIP-based models can handle the uncertainties (Neri et al., 2022; Yi et al., 2021). MIP formulations may also be divided into Mixed-Integer Linear Programming (MILP) and Mixed-Integer Nonlinear Programming (MINLP). The optimal scheduling and operation by stochastic MILP models are frequently studied (Abunima et al., 2022; Dini et al., 2022; Dong et al., 2023; Fusco et al., 2023; Jamalzadeh et al., 2020; Kumar et al., 2023; Thang et al., 2022; Tostado-Véliz et al., 2022; Vergara et al., 2020; Yuan et al., 2020). On the other hand, stochastic MIP-based nonlinear models for the simultaneous optimal design, configuration, and scheduling, are scarce points in the literature (Hafiz et al., 2019; Lee et al., 2022; Vera et al., 2023; Zhou et al., 2013).

This section summarizes the recent studies based on stochastic MIP models for optimal design and scheduling, or only for optimal scheduling of energy grids. Firstly, the examples of studies based on only optimal scheduling and operation are summarized. For example, a novel scenario-based stochastic MINLP formulation is proposed for the optimal flexible-reliable operation of energy hubs for electricity, natural gas, and heating networks, which include renewable energy sources (PV and wind systems), a CHP system, and energy storage systems by incentive-based demand response program. The model uncertainties are dependent on load, energy cost, power generation of renewable-based generators, and network equipment availability (Dini et al., 2022). Moreover, an MINLP-based stochastic scheduling model is suggested to meet electrical, heating, and cooling demands. The uncertainty sources are solar energy, loads, and energy prices (Thang et al., 2022). Furthermore, a linearized form of a stochastic MINLP-based model is used for the optimal operation of islanded microgrids, which comprise an unbalanced three-phase electrical distribution, battery systems and wind turbines to meet the electricity demand (Vergara et al., 2020). The uncertainty sources are stochastic demands and renewable resources. In (Dong et al., 2023), novel forecast-driven mixed robust stochastic and stochastic MINLP-based models are formulated for optimal operation of an isolated microgrid to meet the electricity demand. Wind power and load are uncertainty sources. In (Tostado-Véliz et al., 2022), a MILP-based model is suggested for optimal scheduling of isolated microgrids considering components' failures to meet the electricity demand. The uncertainties sources are weather conditions, wind speed, and demand. In (Budiman et al., 2022), a stochastic MILP-based model is offered for optimal scheduling of a renewable energy-based microgrid with a hybrid energy storage system consisting of a battery and a super-capacitor to meet the electricity demand. The uncertainties sources are wind power, solar power, and load. In (Eghbali et al., 2022), a stochastic MILP-based model is formulated for the optimal scheduling of smart microgrids, which include a wind turbine, a photovoltaic unit, a fuel cell, an electrolyzer, a microturbine, and energy storage units to meet the electricity demand. The uncertainty sources are wind speed, solar irradiance, demand, and electricity price. In (Abunima et al., 2022), a two-stage stochastic MINLP-based model is formulated for optimal operation of a renewable-based microgrid to meet the electricity demand. The uncertainty source is solar energy. In (Fusco et al., 2023), a multi-stage stochastic MILP-based model with binary recourse is proposed for optimizing the day-ahead unit commitment of power plants and virtual power plants to meet the electricity and heat demand. The uncertainty sources are the ancillary services market demands and solar energy. In (Kumar et al., 2023), a tri-level stochastic MINLP-based optimization is used for the optimal operation of grid-connected and autonomous microgrids to eliminate power losses and maximize loadability to meet the electricity demand. The uncertainty sources are solar radiation, wind speed, load, and power loss.

On the other hand, only a few of studies focus on the optimal design, configuration, and scheduling by stochastic MIP, and in particular in nonlinear form (Tatar and Aydin, 2024). For example, a geographic-information-based two-stage stochastic MILP model is formulated for optimal design, location, and scheduling of a MEMG to meet the electricity demand. The uncertainty sources are solar power and demand (Vera et al., 2023). In (Lee et al., 2022), a two-stage stochastic MILP model is proposed for optimal capacity and operation of an electrochemical carbon dioxide reduction system coupled with renewable energy generators and batteries to meet the electricity and methanol demands. The uncertainty sources are weather conditions. In (García-Muñoz et al., 2022), a linearized form of a two-stage stochastic MINLP model is suggested for optimal sizing, allocation, and operation of distributed energy generators coupled with electric vehicles to meet the electricity demand. The uncertainty sources are wind speed, solar radiation, and load.

As seen in the aforementioned studies, there is a significant gap in the literature about stochastic MINLP studies for the optimal equipment selection to simultaneously design and schedule energy grids (Hafiz et al., 2019; Lee et al., 2022; Vera et al., 2023; Zhou et al., 2013). Hence, our study offers a two-stage stochastic MINLP model for the optimal equipment selection to design and plan a MEMG. This way, the best equipment design and sizing decisions that are optimally responding to uncertain scenarios. Moreover, it is a novel study of a two-stage stochastic model that considers uncertainties related to implementing carbon dioxide regulation policies, which are changeable over the years in a country. Another contribution is that this study compares the optimal equipment selections between deterministic and stochastic forms of the MINLP model to see the impact of uncertainties. Since energy grid design at the core of P2G units is scarce in the literature (Carrera and Azzaro-Pantel, 2021; Janke et al., 2022; Pastore et al., 2022), this study also considers the optimal equipment selections to couple with a P2G system under uncertainties as a case study. The uncertainty sources are also more comprehensive than the aforementioned studies (Abunima et al., 2022; Dini et al., 2022; Dong et al., 2023; Fusco et al., 2023; Hafiz et al., 2019; Jamalzadeh et al., 2020; Kumar et al., 2023; Lee et al., 2022; Thang et al., 2022; Tostado-Véliz et al., 2022; Vera et al., 2023; Vergara et al., 2020; Yuan et al., 2020; Zhou et al., 2013).

The major contributions of the proposed work are summarized as:

- Proposing a novel two-stage stochastic MINLP model considering uncertainties related to the implementation of carbon dioxide regulation policies, which are carbon trading and carbon dioxide emission.

- Suggesting the optimal equipment selections, design, and scheduling for a MEMG just in one formulation of a two-stage stochastic MINLP model.

- Offering a MEMG design to couple with a P2G system under uncertainties.

- Including wide variety uncertainty sources, such as carbon dioxide regulation policies (carbon trading and carbon dioxide emission taxing), weather conditions (solar radiation, air temperature, and wind speed), carbon dioxide taxing or trading price, and natural gas price.

The remaining of the paper is organized as follows: Section 3 represents the proposed two-stage stochastic MINLP-based model, candidate equipment modeling and the case studies. Section 4 shows the results and discussion part. Finally, Section 5 summarizes the conclusion points of the study.

## 3. Methods

### 3.1. Description of the Optimization Problem

The carbon trading price is assumed to be the same as the carbon dioxide for fair comparison between carbon trading and emission taxing. The installation location of the MEMG is in Antakya, Türkiye. Utilized seasonal data of hourly wind speed, solar irradiation, and air temperature are obtained from actual measurements (Akulker and Aydin, 2024, 2023). Thirty-two scenarios are created based on the uncertainty sources, as seen in Table 1. The carbon dioxide limit is calculated based on 0.3 tons per generated electricity (Akulker and Aydin, 2023). In addition, the upper and lower limits of carbon dioxide and natural gas prices are determined according to fluctuations (Akulker and Aydin, 2024, 2023).

**Table 1.** Created scenarios.

| Scenarios (w) | Data Type | Cap & Trade | Emission Taxing | Carbon Dioxide Price or Tax ($/ton) | Natural Gas Price ($/m³) |
|---|---|---|---|---|---|
| w1  | Winter | ✓ |   | 50  | 0.86 |
| w2  | Winter | ✓ |   | 100 | 0.86 |
| w3  | Winter |   | ✓ | 50  | 0.86 |
| w4  | Winter |   | ✓ | 100 | 0.86 |
| w5  | Spring | ✓ |   | 50  | 0.86 |
| w6  | Spring | ✓ |   | 100 | 0.86 |
| w7  | Spring |   | ✓ | 50  | 0.86 |
| w8  | Spring |   | ✓ | 100 | 0.86 |
| w9  | Summer | ✓ |   | 50  | 0.86 |
| w10 | Summer | ✓ |   | 100 | 0.86 |
| w11 | Summer |   | ✓ | 50  | 0.86 |
| w12 | Summer |   | ✓ | 100 | 0.86 |
| w13 | Autumn | ✓ |   | 50  | 0.86 |
| w14 | Autumn | ✓ |   | 100 | 0.86 |
| w15 | Autumn |   | ✓ | 50  | 0.86 |
| w16 | Autumn |   | ✓ | 100 | 0.86 |
| w17 | Winter | ✓ |   | 50  | 1.72 |
| w18 | Spring |   | ✓ | 100 | 1.72 |
| w19 | Summer | ✓ |   | 50  | 1.72 |
| w20 | Autumn |   | ✓ | 100 | 1.72 |
| w21 | Winter |   | ✓ | 50  | 1.72 |
| w22 | Spring | ✓ |   | 100 | 1.72 |
| w23 | Summer |   | ✓ | 50  | 1.72 |
| w24 | Autumn | ✓ |   | 100 | 1.72 |
| w25 | Winter |   | ✓ | 100 | 1.72 |
| w26 | Spring | ✓ |   | 50  | 1.72 |
| w27 | Summer |   | ✓ | 100 | 1.72 |
| w28 | Autumn | ✓ |   | 50  | 1.72 |
| w29 | Winter | ✓ |   | 100 | 1.72 |
| w30 | Spring |   | ✓ | 50  | 1.72 |
| w31 | Summer | ✓ |   | 100 | 1.72 |
| w32 | Autumn |   | ✓ | 50  | 1.72 |

Each day is identical over the year. Each day comprises same hourly wind speed, air temperature, and solar radiation for each individual scenario. For example, average winter data are used under carbon trading policy with the carbon dioxide price of 50 $/ton and the natural gas price of 0.86 $/m³ for scenario 1 (*w1*).

The operating powers (power outputs), rated capacities, and carbon dioxide emission of each piece of candidate equipment are required to solve the optimization problem. Table 2 shows all related information for each part of the candidate equipment. However, the model equations are detailed only for the candidate equipment, which has nonlinear constraints and decision variables.

**Table 2.** Some specifications for the candidate equipment (Akulker and Aydin, 2023; Tatar et al., 2022).

| Candidate Equipment | Min. Rated Power (MW) | Max. Rated Power (MW) | Min. Storage Capacity (MWh) | Max. Storage Capacity (MWh) | Installation Cost per Rated Power (M$/MW) | Installation Cost per Storage Capacity ($/MWh) | Maintenance Cost per Rated Power ($/MW) | Maintenance Cost per Storage Capacity ($/MWh) | $CO_2$ Generation (ton/MWh) |
|---|---|---|---|---|---|---|---|---|---|
| WT-1 | 35.4 | 61 | - | - | 2.31 | - | 46,523 | - | 0.025 |
| WT-2 | 49.1 | 108 | - | - | 2.31 | - | 46,523 | - | 0.025 |
| SPA-1 | 16.5 | 16.5 | - | - | 2.39 | - | 28,685 | - | 0.038 |
| SPA-2 | 16.5 | 16.5 | - | - | 2.39 | - | 28,685 | - | 0.038 |
| BCC-1 | 6.6 | 6.6 | - | - | 3.97 | - | 310,777 | - | 0.107 |
| BCC-2 | 11.6 | 11.6 | - | - | 3.31 | - | 235,177 | - | 0.107 |
| BBFB | 2.5 | 50 | - | - | 3.24 | - | 209,354 | - | 0.079 |
| CVT-1 | 10 | 50 | - | - | 0,83 | - | 41,856 | - | Eq. (16) |
| CVT-2 | 10 | 50 | - | - | 0,83 | - | 41,856 | - | Eq. (16) |
| CHP-1 | 20 | 60 | - | - | 1.89 | - | 69,231 | - | 0.840 |
| CHP-2 | 35 | 105 | - | - | 1.903 | - | 86,538 | - | 0.750 |
| CHP-3 | 40 | 125.8 | - | - | 1.9 | - | 85,000 | - | 0.760 |
| IGCC-1 | 10 | 50 | - | - | 2.89 | - | 95,782 | - | 0.319 |
| IGCC-2 | 10 | 50 | - | - | 2.60 | - | 87,438 | - | 0.319 |
| ES | 1 | 50 | 100 | 500 | - | 1.4 | 15,508 | 38,769 | 0 |
| HP | 0 | 50 | - | - | 1.60 | - | 48,000 | - | 0 |
| P2G | 10 | 10 | - | - | 3.258 | - | 114,000 | - | 0 |

## 3.2. Wind Turbine Power Generation Modeling

Each farm contains ten identical wind turbines. The rotor diameter of each wind turbine is the decision variable. The wind turbines in the market, whose diameters are specific, are not preferred as nominees. Instead, in this study, the diameters of wind turbines are optimized for production design. The operational and rated powers of wind turbines depend on their rotor diameters. Hence, the model equations for their operating and rated powers are explained in this section.

Before the operating power (power out-put) calculation, the actual wind speed at the installation area should be found in Eq. (1).

$$v_{wind} = v_{Anemometer} \left(\frac{Z_{WT}}{Z_{Anemometer}}\right)^{\alpha} \tag{1}$$

$v_{Anemometer}$ is the wind speed at the height of the anemometer above the sea level, $Z_{Anemometer}$. $Z_{WT}$ is the height the installation location. $\alpha$ is the friction coefficient (Akulker and Aydin, 2023; Masters, 2004).

The operating power ($p_{WT}$) and the rated power of each wind turbine ($rp_{WT}$) are formulated in Eqs. (2) and (3), respectively.

$$p_{WT} = \begin{cases} 0 & , v_{wind} < v_{cut-in} \\ \frac{1}{2}\rho_{air}\left(\pi\frac{D_{WT}^2}{4}\right)C_p (v_{wind})^3 & , v_{cut-in} \leq v_{wind} < v_{rated} \\ rp_{WT} & , v_{wind} \geq v_{rated} \\ 0 & , v_{wind} \geq v_{cut-out} \end{cases} \quad (2)$$

$$rp_{WT} = \frac{1}{2}\rho_{air}\left(\pi\frac{D_{WT}^2}{4}\right)C_p (v_{rated})^3 \quad (3)$$

$\rho_{air}$ is the standard air density, and $C_p$ is rotor efficiency (Masters, 2004; Wass, 2018).

### 3.3. Solar Panel Power Generation Modeling

Mono-crystalline and polycrystalline PV cells are preferred for SPA-1 and SPA-2, respectively. Tilted module configuration is chosen for design. The tilt angle of each array ($\beta_{SPA}$) is the decision variable for the optimization model, and it is let to be between 20° and 70°.

The actual operating power of a tilted PV cell is calculated in Eq. (4).

$$p_{SPA} = G_{\beta,SPA}\eta_{SPA}A_{SPA}\eta_{SPA,inverter} \quad (4)$$

$G_{\beta,SPA}$ is the hourly total solar radiation on a tilted PV cell, $\eta_{SPA}$ is the cell efficiency, $A_{SPA}$ is the area, and $\eta_{SPA,inverter}$ is the inverter efficiency. $\eta_{SPA}$ is dependent on $G_{\beta,SPA}$. Thus, once $G_{\beta,SPA}$ is determined, $\eta_{SPA}$ is calculated.

$G_{\beta,SPA}$ is determined in Eq. (5).

$$G_{\beta,SPA} = G_{d\beta,SPA} + G_{b\beta} + G_{r\beta,SPA} \quad (5)$$

$G_{d\beta,SPA}$ is hourly diffuse solar radiation on a tilted plane, $G_{b\beta}$ is hourly beam radiation, and $G_{r\beta,SPA}$ is hourly reflected radiation. $G_{d\beta,SPA}$ is dependent on clearness index (*CI*), diffuse fraction (*f*), and diffuse solar radiation on a horizontal plane (*G_d*), which are determined as in Eqs. (6), (7), and (8), respectively (Perera et al., 2012).

$$CI = \frac{G}{H_0} \tag{6}$$

$H_0$ is the extraterrestrial solar radiation.

$$f = \begin{cases} 0.995 - 0.081 CI & , CI < 0.21 \\ 0.724 + 2.738 CI - 8.32 CI^2 + 4.967 CI^3 & , 0.21 \le CI \le 0.76 \\ 0.18 & , 0.76 < CI \end{cases} \tag{7}$$

$$G_d = f \cdot G \tag{8}$$

$G_{d\beta,SPA}$ is calculated by the Klucher model (Klucher, 1979):

$$F = 1 - f^2 \tag{9}$$

$$G_{d\beta,SPA} = G_d[0.5(1 + \cos(\beta_{SPA}/2))][1 + F\sin^3(\beta_{SPA}/2)] \cdot [1 + F\cos^2(\theta)\sin^3(\theta_z)] \tag{10}$$

$\theta$ is the declination angle presenting the angle between the sun lights and equator plane. $\theta_z$ is the zenith angle (Hailu and Fung, 2019).

Beam ($G_{b\beta}$) and reflected radiation ($G_{r\beta,SPA}$) are calculated as in Eqs. (11) and (12), respectively (Klucher, 1979; Perera et al., 2012).

$$G_{b\beta} = (G - G_d)\cos(\theta)/\cos(\theta_z) \tag{11}$$

$$G_{r\beta,SPA} = \rho G(1 - \cos(\beta_{SPA}))/2 \tag{12}$$

$\rho$ is the Albedo ratio.

The hourly PV cell efficiency ($\eta_{SPA}$) (%) is calculated in Eq. (13).

$$\eta_{SPA} = p\_spa \left[ q\_spa \frac{G_{\beta,SPA}}{G_{\beta 0}} + \left( \frac{G_{\beta,SPA}}{G_{\beta 0}} \right)^{m\_spa} \right] \left[ 1 + r\_spa \frac{\theta_{cell\_spa}}{\theta_{cell,0}} + s\_spa \frac{AM}{AM_0} \right. \quad (13)$$
$$\left. + \left( \frac{AM}{AM_0} \right)^{u\_spa} \right]$$

*AM* is the air mass value (Kasten and Young, 1989), and $\theta_{cell\_spa}$ is the cell temperature, calculated by Eq. (14). $G_{\beta 0}$, $\theta_{cell,0}$, and $AM_0$ are 1000 W/m², 25°C, and 1.5, respectively.

$$\theta_{cell\_spa} = \theta_a + h\_spa \, G_{\beta,SPA} \quad (14)$$

$\theta_a$ is the ambient temperature. *p_spa, q_ spa, m_ spa, r_ spa, s_ spa, u_ spa,* and *h_ spa* are the constant coefficients obtained from (Akulker and Aydin, 2023; Perera et al., 2012).

### 3.4. Conventional generator modeling

The hourly operating cost ($OC_{cvt}$) function for conventional generators is nonlinearly dependent on their operating powers ($p_{cvt}$) as follows (Zhang et al., 2013):

$$OC_{cvt} = c\_cvt(p_{cvt}^2) + b\_cvt(p_{cvt}) + a\_cvt \quad (15)$$
$$+ \left| d\_cvt \left( sin \left( e_{cvt}(p_{cvt}^{min} - p_{cvt}) \right) \right) \right|$$

*a_cvt, b_cvt, c_cvt, d_cvt,* and *e_cvt* are the constant coefficients obtained from (Akulker and Aydin, 2023). $p_{cvt}^{min}$ is the minimum operating power.

The hourly carbon dioxide emission from the conventional generators is calculated in Eq. (16).

$$CO2_{cvt} \left( \frac{ton}{hour} \right) = ef(h\_cvt(p_{cvt})^2 + g\_cvt(p_{cvt}) + f\_cvt) \quad (16)$$

*ef, h_cvt, g_cvt,* and *f_cvt* are the constant coefficients of the hourly carbon dioxide emission function, which are taken from (Akulker and Aydin, 2023; Zhang et al., 2013).

## 3.5. Combined heat and power (CHP) generator modeling

All candidate CHP units have tetragonal feasible operational regions, as listed in Appendix Table A1. Fig. 1 shows the feasible operational region of CHP-2.

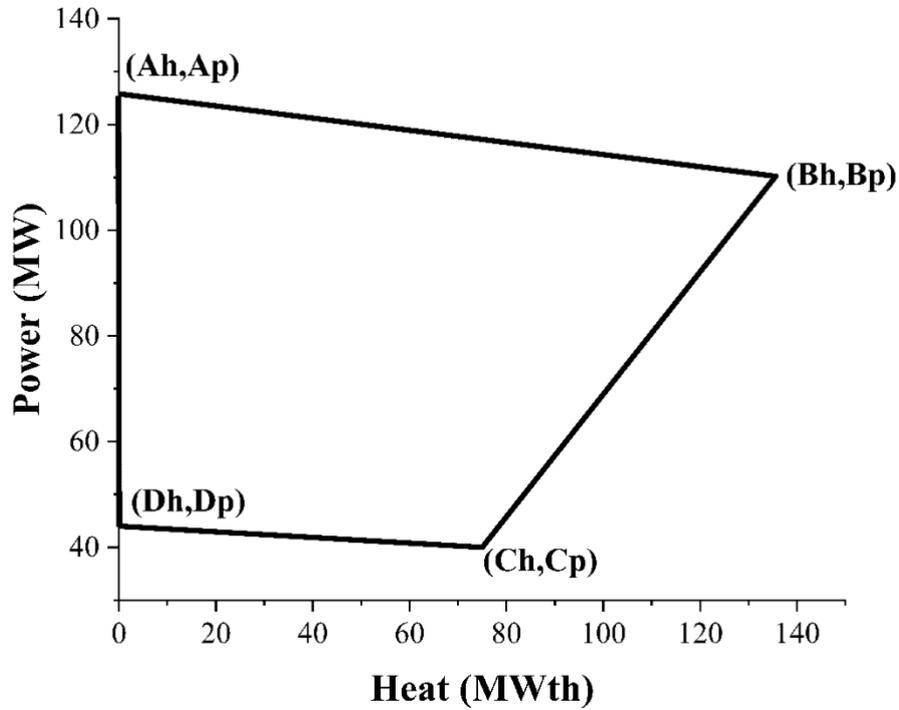

**Fig. 1.** The feasible operational region for CHP-2.

Feasible operation region constraints are integrated in to the optimization model according to Fig. 1, as given in Eqs. (17)-(19) (Ko and Kim, 2019).

$$\left(p_{chp} - 'D_p'\right) \geq \left(h_{chp} - 'D_h'\right) \cdot \left(\frac{('D_p' - 'C_p')}{('D_h' - 'C_h')}\right) \tag{17}$$

$$\left(p_{chp} - 'A_p'\right) \leq \left(h_{chp} - 'D_h'\right) \cdot \left(\frac{('A_p' - B_p')}{('A_h' - 'B_h')}\right) \tag{18}$$

$$(p_{chp} - 'B_p') \geq (h_{chp} - 'B_h') \cdot \left(\frac{('B_p' - 'C_p')}{('B_h' - 'C_h')}\right) \quad (19)$$

$p_{chp}$ and $h_{chp}$ are the operational power and heat of CHP units, respectively. $D_h$ and $D_p$ are the coordinates of the corner *D*. The similar notations are given for *A, B, C*.

The operational (fuel consumption) cost function of CHP generators is calculated as in Eq. (21).

$$OC_{chp}\left(\frac{\$}{\text{hour}}\right) = kk_{chp} + ll_{chp}(p_{chp}) + ii_{chp}(p_{chp})^2 + jj_{chp}h_{chp} + yy_{chp}(h_{chp})^2 + zz_{chp}h_{chp}p_{chp} \quad (20)$$

*kk_chp, ll_chp, ii_chp, jj_chp, yy_chp*, and *zz_chp* are the operational constant cost coefficients, as shown in Appendix Table A2.

### 3.6. Power-to-Gas (P2G) system modeling

The P2G system includes a solid oxide electrolyzer cell (SOEC) to produce hydrogen, a carbon capture system (CCS) to extract carbon dioxide, a methanation unit including a fixed-bed reactor with nickel-alumina catalyst, a separator, and an upgrading unit for methane to SNG (Gorre et al., 2019).

The exothermic Sabatier reaction occurs in the methanation unit, as given in Eq. (21).

$$CO_2 + 4H_2 \leftrightarrow CH_4 + 2H_2O, \quad \Delta H° = -165 \; kJ/mol \quad (21)$$

Thermodynamic equilibrium is the constraint of the reaction as shown in Eq. (22).

$$PP_{CO2}(PP_{H2})^4 K_{eq}(T^{in}) \leq (PP_{H2O})^2 PP_{CH4} \quad (22)$$

$PP_{CO2}$, $PP_{H2}$, $PP_{H2O}$, and $PP_{CH4}$ are the partial pressures. $K_{eq}$ is the equilibrium constant (Davis and Martín, 2014), and $T^{in}$ is the inner temperature of the reactor.

The partial pressures of the species are calculated as in Eqs. (23) and (24).

$$PP_j = PP^{in} \frac{N_j}{\sum_j N_j} \tag{23}$$

$$N_j = N_j^{in} + \xi \nu_j \quad \text{for } j=\{CO_2, H_2, H_2O, CH_4\} \tag{24}$$

$PP_j$ is the partial pressure for the species $j$. $N_j$ is the molar output flow rate of the species $j$. $\nu_j$ is the stoichiometric coefficient of the species $j$. $N_j^{in}$ is the molar inlet flow rate of the species $j$, and $PP^{in}$ is the inner pressure of the reactor. $\xi$ is the extent of the reaction (Uebbing et al., 2020).

Fig. 2 shows the superstructure of the MEMG.

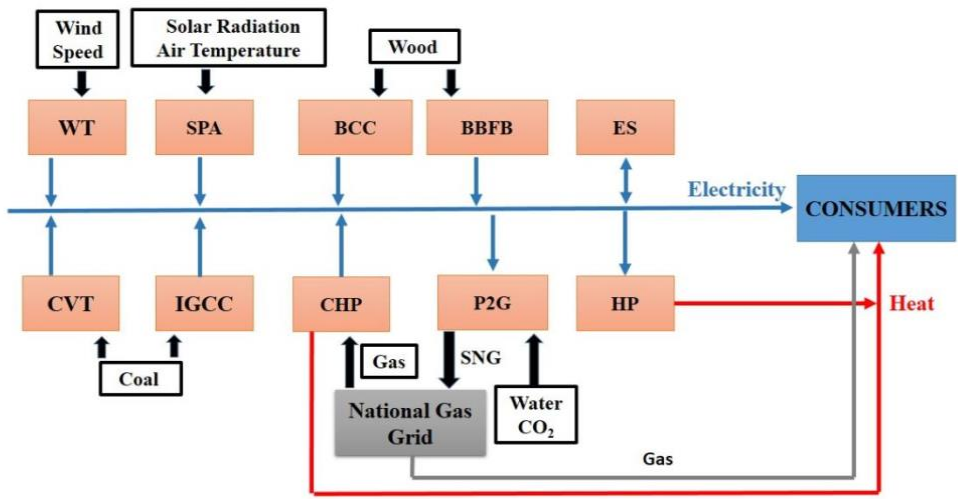

**Fig. 2.** The superstructure of the MEMG.

## 3.7. The Two-Stage Stochastic MINLP Optimization Model

The optimal design and operation of MEMGs by handling uncertainties are formulated as a two-stage stochastic programming problem. The first stage is related to the design variables responding to all realizations of uncertainty, whereas the second-stage is an operational-level problem also called as recourse problem, following the realization of first-stage design-level decision variables. All design variables are here-and-now variables, whereas all the operation variables are wait-and-see variables. The objective function is to minimize the total annualized cost (TAC) of the installation and operation of the MEMG under uncertainties. The two-stage stochastic formulation of the objective function is written in Eqs. (25)-(27):

**First-stage:**

$$\min \left((CRF)\, C^{Capital}(i, rp, cap) + C^{Maintenance}(i, rp, cap) + \mathbb{E}\left[C^{O,CT,ET,R}(i, rp, cap, \varepsilon_f)\right]\right) \tag{25}$$

$$\text{s.t. } \emptyset(rp, cap) \leq 0$$

$$CRF = \left(\frac{int(int+1)^{LFT}}{(int+1)^{LFT}-1}\right) \tag{26}$$

**Second-stage:**

$$C^{O,CT,ET,R}(i, rp, cap, \varepsilon_f) = \min(C^{Operational}(i, rp, cap, p, kc, \varepsilon_f) + C^{CT}(i, rp, cap, p, kc, \varepsilon_f) + C^{ET}(i, rp, cap, p, kc, \varepsilon_f) - C^{Revenue}(i, rp, cap, p, kc, \varepsilon_f)) \tag{27}$$

$$\text{s.t. } \mathsf{Y}(rp, cap, p, \varepsilon_f) \leq 0$$

$$\emptyset(rp, cap) \leq 0$$

"*i*" is the binary variable to decide the candidate equipment to be installed. *rp* and *cap* are the rated power of each piece of candidate equipment and the storage capacity of each candidate storage unit, respectively. $\mathbb{E}$ is the expectation of stochastic costs and revenue ($C^{O,CT,ET,R}$). In this study, $\varepsilon_f$ is the realization of uncertainties related to seasonal data (hourly air temperature, solar radiation, and wind speed), carbon dioxide regulation policies (carbon trading and emission taxing), carbon dioxide price or tax, and natural gas price. As seen in Eq. (26), $C^{Capital}$ (the total installation cost) and $C^{Maintenance}$ (the total maintenance cost) are related to the first-stage design-level variables. They are dependent only on the equipment selections (*i*), the rated powers of the candidate equipment (*rp*), and the storage capacity of the candidate batteries (*cap*). In addition, they are subject to fixed and structural constraints ($\emptyset(rp, cap) \leq 0$). Installation decisions must be made before the values of uncertain parameters are observed. In other words, stochastic programming selects equipment for installation to withstand all uncertainties defined in the optimization model. *CRF* is the capital recovery factor, as calculated in Eq. (26). *int* is interest rate (15%) and *LFT* is the lifetime of each candidate equipment (20 years).

Since the uncertainties affect the operating power (*p*) of each candidate piece of equipment, $C^{Operational}$ (operational cost), $C^{CT}$ (carbon trading cost or revenue), $C^{ET}$ (emission taxing cost), and $C^{Revenue}$ (the revenue from the sale of synthetic natural gas (SNG)) are also affected. "*kc*" is the binary variable to schedule the corresponding equipment, in other words to decide which candidate equipment to switch on or off, and $\gamma(rp,b,p,\varepsilon_f)$ is the general description of the constraints related to uncertain parameters. The second-stage variables represent recourse actions, also called recourse variables (Barbarosoğlu and Arda, 2004; Yang et al., 2017).

The two-stage stochastic formulation is converted to its deterministic equivalent to solve the problem with the finite number of scenarios. The objective function of the deterministic equivalent of the optimization problem is given as follows:

$$\min (CRF) C^{Capital}(i, rp, cap) + C^{Maintenance}(i, rp, cap)$$
$$+ \sum_w \pi_w (C^{Operational}(i, rp, cap, p, k, \varepsilon_{f,w})$$
$$+ C^{CT}(i, rp, cap, p, k, \varepsilon_{f,w}) + C^{ET}(i, rp, cap, p, k, \varepsilon_{f,w})$$
$$- C^{Revenue}(i, rp, cap, p, k, \varepsilon_{f,w})) \quad (29)$$

$\pi_w$ is the probability of each scenario ($w$). $\varepsilon_{f,w}$ is the realization of the uncertainties for each scenario, $w$. The probability of each scenario is assumed to be equal. Thus, $\pi_w$ is accepted as 1/32 for the calculations.

The rated power of each piece of candidate equipment, $rp(c)$, and the storage capacity of each battery, $cap(c)$ are continuous decision variables. Their limits are described in Eqs. (30) and (31).

$$rp_c^{min} i(c) \leq rp(c) \leq rp_c^{max} i(c) \quad for \ \forall c \in C \ and \ i(c) \in \{0,1\}. \quad (30)$$

$$cap_c^{min} i(c) \leq cap(c) \leq cap_c^{max} i(c) \quad (31)$$

$for \ \forall c \in ST \in C \ and \ i(c) \in \{0,1\}.$

$rp_c^{min}$, $rp_c^{max}$, $cap_c^{min}$, and $cap_c^{max}$ are the lower and upper limits for the rated powers and the battery capacity, respectively, as shown in Table 2.

$kc_{c,w,t}$ decides whether the piece of candidate equipment, $c$, is operated at time interval $t$ in scenario $w$ as in Eqs. (32) and (33).

$$p_{ct}^{min}(rp(c) - (1 - kc_{c,w,t})rp_c^{max}) \leq p_{c,w,t} \leq p_{ct}^{max}rp(c) \tag{32}$$

$$0 \leq p_{c,w,t} \leq p_{ct}^{max}rp_c^{max}kc_{c,w,t} \tag{33}$$

$for\ \forall t \in T, \forall c \in G \in C, \forall w \in W, and\ kc_{c,w,t} \in \{0,1\}.$

For storage units, the charging ($pch_{c,w,t}$) and discharging ($pdch_{c,w,t}$) powers should be between 0% ($q^{min}$) and 100% ($q^{max}$) of their rated powers, as seen in Eqs. (34) and (35).

$$q^{min}rp(c) \leq pch_{c,w,t} \leq q^{max}rp(c) \quad for\ \forall c \in ST \in C,\ \forall t \in T, and\ \forall w \in W. \tag{34}$$

$$q^{min}rp(c) \leq pdch_{c,w,t} \leq q^{max}rp(c) \quad for\ \forall c \in ST \in C,\ \forall t \in T, and\ \forall w \in W. \tag{35}$$

In addition, Eqs. (36) and (37) indicate the constraints that the charging and discharging of the storage unit cannot operate at the same time interval.

$$pch_{c,w,t} \leq q^{max}rp_c^{max}kc_{c,w,t} \tag{36}$$

$$pdch_{,sn,t} \leq q^{max}rp_c^{max}(1 - kc_{c,w,t}) \tag{37}$$

$for\ \forall t \in T, \forall c \in ST \in C, \forall w \in W, and\ kc_{c,w,t} \in \{0,1\}.$

Moreover, energy stored in the storage unit at any time interval is described by $soc_{c,w,t}$ is given in Eq. (38) by regarding safety and equipment lifetime.

$$(20\%)cap(c) \leq soc_{c,w,t} \leq (80\%)cap(c) \quad for \quad \forall c \in ST \in C, \forall t \in \mathcal{T}, \text{and} \quad (38)$$
$$\forall w \in W.$$

Furthermore, Eq. (39) shows the capacity limits, the initial and final conditions of the battery.

$$soc_{c,w,t} = \begin{cases} soc_{c,w,0} + (pch_{c,w,t} - pdch_{c,w,t}), & if\ t = t(1) \\ soc_{c,w,'t-1'} + (pch_{c,w,t} - pdch_{c,w,t}), & otherwise \end{cases} \quad (39)$$

$$for \quad \forall c \in ST \in C, \forall t \in \mathcal{T}, \text{and } \forall w \in W.$$

Eq. 40 indicates power and material balance of the MEMG. The generation and consumption terms include linear and nonlinear equations. Therefore, equipment $C$ set is split into two subsets, which are *CL* and *CNL,* to express the whole power and material balance just in a single formulation. *CL* comprises the candidate equipment owing to only linear constraints dependent on decision variables. In contrast, C*NL* includes the candidate equipment having both linear and nonlinear constraints. WTF-1, WT-2, SPA-1, SPA-2, CVT-1, CVT-2, CHP-1, CHP-2, CHP-3 and P2G are the members of *CNL,* whilst the others are of *CL*.

$$\sum_{cnl} gen_{cnl,w,t,res} + \sum_{cl} gen_{cl,w,t,res} + u_{w,t,res} \quad (40)$$
$$= \sum_{cnl} con_{cnl,w,t,res} + \sum_{cl} con_{cl,w,t,res} + yx_{w,t,res} + spin_{w,t,res}$$
$$+ demand_{w,t,res}$$

$for\ \forall t \in \mathcal{T},\ \forall cl \in CL \in C,\ \forall cnl \in CNL \in C,\ \forall w \in W, and\ \forall res \in RES.$

$con_{cnl,w,t,'coal'}$ is considered zero for the conventional power generators because the operational costs of the conventional generators are determined by a highly nonlinear equation, as shown in Eq. (15). For the same reason, $con_{cnl,w,t,'gas'}$ is taken as zero for the CHP units.

The final constraint is related to the maximum allowable number for equipment installation, as in Eq. (41).

$$\sum_{c \in C} i(c) \leq 9 \tag{41}$$

$C^{Capital}$ is calculated in Eq. (42).

$$C^{Capital} = \sum_{c \in C} (\psi_c^0 rp(c) + \gamma_c^0 a(c)) + \sum_{c \in ST} \omega_c^0 b(c) \tag{42}$$

$\psi_c^0$ and $\omega_c^0$ are the installation coefficients of candidate equipment and candidate storage units, respectively. $\gamma_c^0$ is the fixed cost for each piece of candidate equipment. $C^{Maintenance}$ is calculated in Eq. (43).

$$C^{Maintenance} = \sum_{c \in C} \left( \psi_c^k rp(c) + \gamma_c^k a(c) \right) + \sum_{c \in ST} \omega_c^k b(c) \tag{43}$$

$\psi_c^k, \gamma_c^k$, and $\omega_c^k$ are maintenance cost coefficients. All of the cost coefficients are given in Table 2.

$C^{Operational}$ is calculated as in Eq. (44).

$$C^{Operational}(a, rp, b, p, k, \varepsilon_{f,w})$$

$$= 365 \left( \sum_t \sum_{res} u_{w,t,res} \, Prc(res) \right. \tag{44}$$

$$\left. + \sum_t \sum_{chp} OC_{chp,w,t} + \sum_t \sum_{cvt} OC_{cvt,w,t} + \sum_t con_{,PTG',w,t,'CO2'} \, CC \right)$$

$$for \; \forall t \in \mathcal{T}, \; \forall w \in W, \; \forall res \in RES, \; \forall cvt \in CVT, \; and \; \forall chp \in CHP.$$

$Prc(res)$ is the price of input resource *res*. The operational costs (fuel consumption costs) are separately calculated from the costs related to buying other resources, as shown in the first term in Eq. (44), for the reason mentioned in the material and energy balance part. The last term of Eq. (44) indicates the variable operational cost for carbon capturing in CCS.

$C^{CT}(i, rp, cap, p, k, \varepsilon_{f,w})$ is the cost of buying the carbon dioxide limit or the revenue from selling it, which is calculated in Eq. (45). In carbon trading, the rulers specify the carbon dioxide caps for energy producers. For example, when a microgrid exceeds the carbon dioxide cap certified by the rulers, an excess carbon dioxide limit is allowed to buy. On the other hand, when the hourly carbon dioxide emission is lower than the limit, the remaining carbon dioxide limit is permitted to sell (TheWorldBank, 2024.). $L_{CO2,w,t}$ is the hourly carbon dioxide cap in tons of $CO_2$ in scenario w at time interval *t*. The excess carbon dioxide ($yx_{w,t,'CO2'}$) at interval *t* in scenario w is obtained from the material balance given in Eq. (40).

$$C^{CT}(a, rp, cap, p, k, \varepsilon_{f,w}) = \sum_t 365 (yx_{w,t,'CO2'} - L_{CO2,w,t}) Prc(CO_2) \tag{45}$$

$$for \quad \forall t \in \mathcal{T} \; and \; \forall w \in W.$$

$Prc(CO_2)$ is the carbon dioxide trading price per ton of $CO_2$. According to the Paris Agreement based on 0.3 tons of $CO_2$ per operating power of each generator in MW, $L_{CO2,w,t}$ is calculated as in Eq. (46).

$$L_{CO2,w,t} = 0.3 \frac{tons\ of\ CO_2}{MW} \sum_c p_{c,w,t} \tag{46}$$

for $\forall c \in (G \cup WTF \cup SPA), \forall t \in T$, and $\forall w \in W$.

$C^{ET}(i, rp, cap, p, k, \varepsilon_{f,w})$ is the total emission cost for carbon dioxide, as determined in Eq. (47).

$$C^{ET}(a, rp, b, p, k, \varepsilon_{f,w}) = \sum_t 365 (yx_{w,t,'CO2'}) Prc(CO_2) \tag{47}$$

Please notice that carbon dioxide tax and trading price are considered the same.

$C^{Revenue}(i, rp, cap, p, k, \varepsilon_{f,w})$ is the revenue from the sale of the generated SNG to the national gas grid, which is calculated in Eq. (48).

$$C^{Revenue}(i, rp, cap, p, k, \varepsilon_{f,w}) = \sum_t 365\ Prc(SNG)(gen_{'P2G',w,t,'SNG'}) \tag{48}$$

$gen_{'PTG',w,t,'SNG'}$ is the generation of SNG by P2G at interval $t$ in scenario w. $Prc(SNG)$ is the natural gas price.

### 3.8. Case Studies

Three case studies are investigated in this study. The aim of the first case (Case-1) is to select the optimal equipment from all the candidate equipment to meet the electricity and heat demands.

The purpose of the second one (Case-2) is to design and schedule a P2G-integrated MEMG by selecting the optimal equipment to meet the demands of electricity, heat, and natural gas. In this case, the sale of SNG to the national gas is assumed to be mandatory. Hence, hourly SNG demand for all scenarios is considered as 0.250 tons, but higher SNG sale than 0.250 tons is allowed.

For Case-1 and Case-2, the optimal equipment selections and the total annualized costs of the stochastic and deterministic models are compared. The sixteen scenarios are based on carbon trading, while the other sixteen are based on carbon dioxide taxing, as noticed in Table 1. The optimal equipment selections by two-stage stochastic MINLP models with the sixteen scenarios, including only carbon trading, the sixteen scenarios, including only carbon dioxide taxing, and thirty-two scenarios, including both carbon regulation policies. On the other hand, the deterministic models are created separately based on carbon trading and carbon dioxide taxing. The hourly data used in the deterministic models are obtained by averaging all weather data: wind speed, air temperature, and solar radiation. In addition, the average prices are used for natural gas and carbon dioxide tax or trade.

The aim of the third one (Case-3) is to select the optimal equipment to work with sustainable energy generators: wind turbine farms and solar panel arrays. Thirty-two scenarios are included in the two-stage stochastic MINLP model (See Table 1). In this case, the installation of WT-1, WT-2, SPA-1, and SPA-2 is assumed to be obligatory to meet the electricity and heat demands. The another aim of Case-3 is to determine the optimal rotor diameters of wind turbines and the optimal tilt angles of solar panels, which are nonlinear decision variables, when they are forced to be installed. Fig. 3 shows the all details of the case studies.

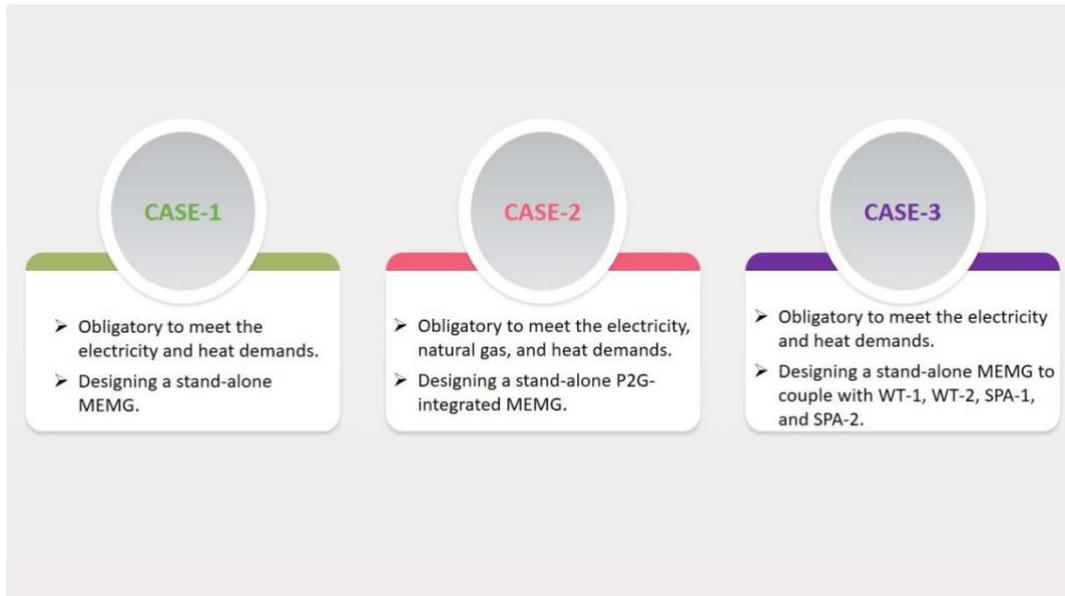

**Fig. 3.** Details of the case studies.

DICOPT, a local solver for non-convex MINLP-based problems, is used in GAMS. DICOPT solves MINLP problems by dividing them into MILP and NLP sub-problems (GAMS, 2024). MILP and NLP solvers are set to CPLEX and CONOPT, respectively. The absolute and relative optimality criteria are set to $10^{-3}$.

## 4. Results and Discussion

### 4.1. The results for Case-1

The optimal equipment selections and the total annualized costs for all the trials in Case-1 are shown in Table 3.

**Table 3.** The optimal equipment selections and the total annualized costs for Case-1.

| Optimization model type | The optimal equipment selections with their rated powers (MW) | TAC (Million $) |
| --- | --- | --- |
| Two-stage stochastic with the sixteen scenarios including only carbon trading | IGCC-1 (50)+ IGCC-2 (50)+ CHP-2 (105)+BBFB (7.711)+ BCC-1 (6.6)+ BCC-2 (11.6) | 278.79 |
| Two-stage stochastic with the sixteen scenarios including only carbon dioxide taxing | IGCC-1 (50)+ IGCC-2 (50)+ CHP-2 (105)+BBFB (7.711)+ BCC-1 (6.6)+ BCC-2 (11.6) | 318.28 |
| Deterministic model including only carbon trading | IGCC-1 (50)+ IGCC-2 (50)+ CHP-2 (105)+WT-2 (83) | 379.05 |
| Deterministic model including only carbon dioxide taxing | IGCC-1 (50)+ IGCC-2 (50)+ CHP-2 (105)+WT-2 (83) | 418.54 |
| Two-stage stochastic with the thirty-two scenarios including both carbon trading and carbon dioxide taxing | IGCC-1 (50)+ IGCC-2 (50)+ CHP-2 (105)+BBFB (7.711)+ BCC-1 (6.6)+ BCC-2 (11.6)+ | 298.54 |

For all the optimization models in Case-1, only CHP-2 is selected to generate heat to meet the whole heat demand because of its double benefit to generate electricity and heat, simultaneously. The heat hump as heat source is not preferred. Heat pumps use excess electricity of the grids to generate heat. Using excess electricity of the MEMG to produce heat is decided not to be profitable by the optimization models in this study. Fig. 4 shows heat generation plan for each scenario.

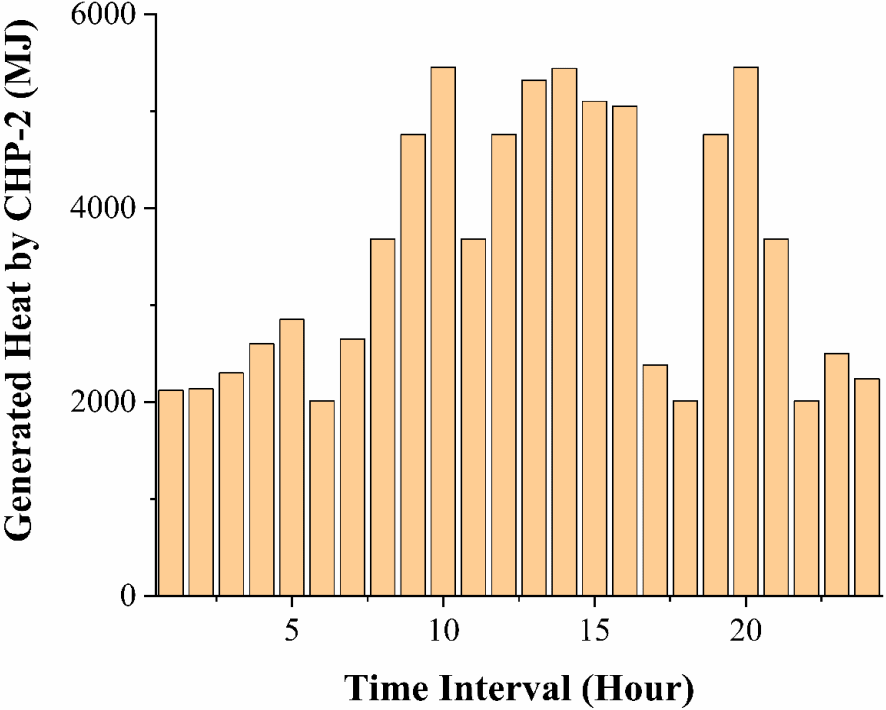

**Fig. 4.** Generated heat by CHP-2 for each scenario in Case-1.

The optimal equipment selections by the two-stage stochastic model, with the sixteen scenarios including only carbon trading system, the sixteen scenarios including only carbon dioxide taxing, and the thirty-two scenarios are the same. IGCCs, BCCs, BBFB, and CHP-2 are selected. Wind turbine farms and solar panel arrays are not chosen, although their carbon dioxide emissions are much lower than the other generators. The intermittency of solar and wind power prevents the continuous energy production over a day. Therefore, the optimization model selects the generators that can maintain energy production in all weather conditions. Although conventional generators can produce energy independently in case of weather, they are not chosen by model. Instead, IGCCs are selected. The reason is that their carbon dioxide emission is much lower than the conventional generators, which also use coal, thanks to the gasification of coal.

On the other hand, the deterministic models, including only carbon dioxide taxing and only carbon trading, choose WT-2 instead of BCCs and BBFB. BCCs and BBFB consume waste wood to generate electricity. Their carbon dioxide emission is much lower than the CHP units and the conventional generators. Since the uncertainty of wind speed is handled in the two-stage stochastic models, this situation is neglected by using the constant average wind speed over the year in the deterministic models; the two-stage stochastic models prefer the generators using more stable energy sources instead of wind turbine farms.

The P2G system is not chosen by any of the optimization models in Case-1 despite implementing strict carbon dioxide regulation policies. Generation of SNG is not profitable for all the trials even though it reduces carbon dioxide emissions.

The TAC values are much lower in the stochastic models than in the deterministic ones. For example, in the stochastic model, the TAC is 278.79 Million $ with the sixteen scenarios including only the carbon trading system, while it is 379.05 Million $ in its deterministic model. The total installation costs of the stochastic and deterministic models are 563.9 Million $ and 666.06 Million $, respectively. The installation of WT-2 instead of BCCs and BBFB increases the total installation cost by nearly 102 Million $. The total operational costs of the stochastic and deterministic models are 138.56 Million $ and 230 Million $, respectively. Operation of WT-2 instead of BCCs and BBFB is expected to lead to highly negligible operational costs. Yet, the probabilistic nature of the stochastic model decreases the total operating cost by nearly 91 Million $. This outcome also explains the difference in the TAC values between the stochastic and deterministic models, including only carbon dioxide emission taxing. The TAC value increases from 318.28 Million $ in the stochastic model to 418.54 Million $ in the deterministic model, as seen in Table 3.

The number of scenarios included in the stochastic model affects the optimal TAC value. Please notice that the TAC, 298.54 Million $, in the stochastic model with thirty-two scenarios is lower than that in the stochastic model with the sixteen scenarios, including only carbon dioxide taxing, 318.28 Million $. However, it is higher than the TAC in the stochastic model with the sixteen scenarios, including only carbon trading system, 278.79 Million $, as seen in Table 3. The installed equipment and their rated powers are the same for these trials. The total operational costs are affected by the number of scenarios for these trials.

Similar to stochastic models, the optimal equipment selections and their rated powers are the same for the deterministic models, including only carbon trading and only carbon dioxide tax. As another significant conclusion, the equipment selections are not affected by the type of carbon regulation policy in the stochastic and deterministic models in this study.

### 4.2. The results for Case-2

The optimal equipment selections and the total annualized costs for all the trials in Case-2 are shown in Table 4.

**Table 4.** The optimal equipment selections and the total annualized costs for Case-2.

| Optimization model type | The optimal equipment selections with their rated powers (MW) | TAC (Million $) |
|---|---|---|
| Two-stage stochastic with the sixteen scenarios including only carbon trading | IGCC-1 (50)+ IGCC-2 (50)+ CHP-2 (105)+BBFB (17.711)+ BCC-1 (6.6)+ BCC-2 (11.6) | 300.23 |
| Two-stage stochastic with the sixteen scenarios including only emission taxing | IGCC-1 (50)+ IGCC-2 (50)+ CHP-2 (105)+BBFB (17.711)+ BCC-1 (6.6)+ BCC-2 (11.6) | 341.69 |
| Deterministic model including only carbon trading | IGCC-1 (50)+ IGCC-2 (50)+ CHP-2 (105)+ WT-1 (35.4) + WT-2 (108.4) | 396.10 |
| Deterministic model including only emission taxing | IGCC-1 (50)+ IGCC-2 (50)+ CHP-2 (105)+ WT-1 (35.4) + WT-2 (108.4) | 437.57 |

Only CHP-2 is selected to generate heat. WT-1 and WT-2 are chosen in the deterministic models instead of BCCs and BBFB. The optimal equipment selections and rated powers of the selected equipment are the same for the deterministic model including only carbon trading and only carbon dioxide taxing.

On the other hand, the optimal equipment selections and their rated powers are the same for the stochastic model including only carbon trading and only carbon dioxide taxing. They both choose BBCs and BBFB instead of wind turbine farms for the same reason explained in Case-1. The optimal equipment selections are highly similar to the selections in the stochastic models in Case-1. The unique difference is in the rated power of BBFB. The P2G-integrated MEMG needs additional electricity. The rated power of the SOEC in the P2G system is 10 MW, so 10 MW must be hourly generated for production of SNG, whose hourly sale to the national gas grid is mandatory to be 0.250 tons. Therefore, the stochastic optimization models in Case-2 choose a higher rated power for BBFB.

Installation of the wind turbine farms instead of BBCs and BBFB increases the TAC in the deterministic models. The TAC in the stochastic model with the sixteen scenarios including only carbon trading system, 300.23 Million $, is much lower than that in its deterministic model, 396.10 Million $.

On the other hand, the TAC in the stochastic model with the sixteen scenarios including only carbon dioxide taxing, 341.69 Million $, is much lower than that in its deterministic model, 437.57 Million $. In addition, the total installation cost is 839.06 Million $ in the deterministic model while it is 628.9 Million $ in the stochastic model. Although the installation cost per rated power for wind turbines is lower than BBCs and BBFB, their rated powers are much higher than the total rated powers of BBCs and BBFB. Please notice that the rated powers of WT-1 and WT-2 are determined as 35.406 MW and 108.383 MW, respectively, whereas the rated powers of BBFB, BCC-1 and BCC-2 are 17.711 MW, 6.6 MW, and 11.6 MW, respectively, as seen in Table 4.

The hourly natural gas (SNG) demand is 0.250 tons. Since the excess sale of SNG is determined to be more profitable for each scenario by the all optimization models in Case-2, 0.492 tons of SNG/hour which is the maximum amount according to the thermodynamic equilibrium, as shown in Eq. (22), are generated. DICOPT is not able to converge within the specified maximum iteration time for the stochastic model with the thirty-two scenarios in this case.

### 4.3. The results for Case-3

The two-stage stochastic model with thirty-two scenarios selects IGCC-1, IGCC-2, BCC-2, CHP-2, whose rated powers are 50 MW, 50 MW, 11.6 MW, and 105 MW, respectively, to couple with sustainable generators: WT-1, WT-2, SPA-1, and SPA-2. The optimal rated powers of WT-1 and WT-2 are found as 35.4 MW and 63.2 MW, respectively. In addition, the optimal rotor diameters for WT-1 and WT-2 are calculated as 76 meter and 114.5 meter, respectively. WT-1 and WT-2 are not installed at their maximum rated powers even though they are forced to be installed.

Since the upper and lower limits of their rated powers are 16.5 MW, the optimal rated powers of SPA-1 and SPA-2 are 16.5 MW when they are forced to be installed. The optimal tilt angle is the continuous decision variable of the optimization model, which directly affects panel efficiencies and operating powers of solar panels. The optimal tilt angle is calculated as 70° for SPA-1 and SPA-2.

The TAC is found as 303.04 Million $. When sustainable generators are not forced to be installed, the same results with the stochastic model with the thirty-two scenarios in Case-1 are obtained. In that trial, the TAC is determined as 298.54 Million $. The difference is roughly 5 Million $, but in terms of sustainable and green energy production, this amount may be compensated to reach net-zero carbon emission target.

Only CHP-2 is selected to generate heat. Hourly electricity generation plan for each scenario in each case could be plotted. However, for brevity, only hourly electricity generation plan for scenario *w31* is plotted in Fig. 5.

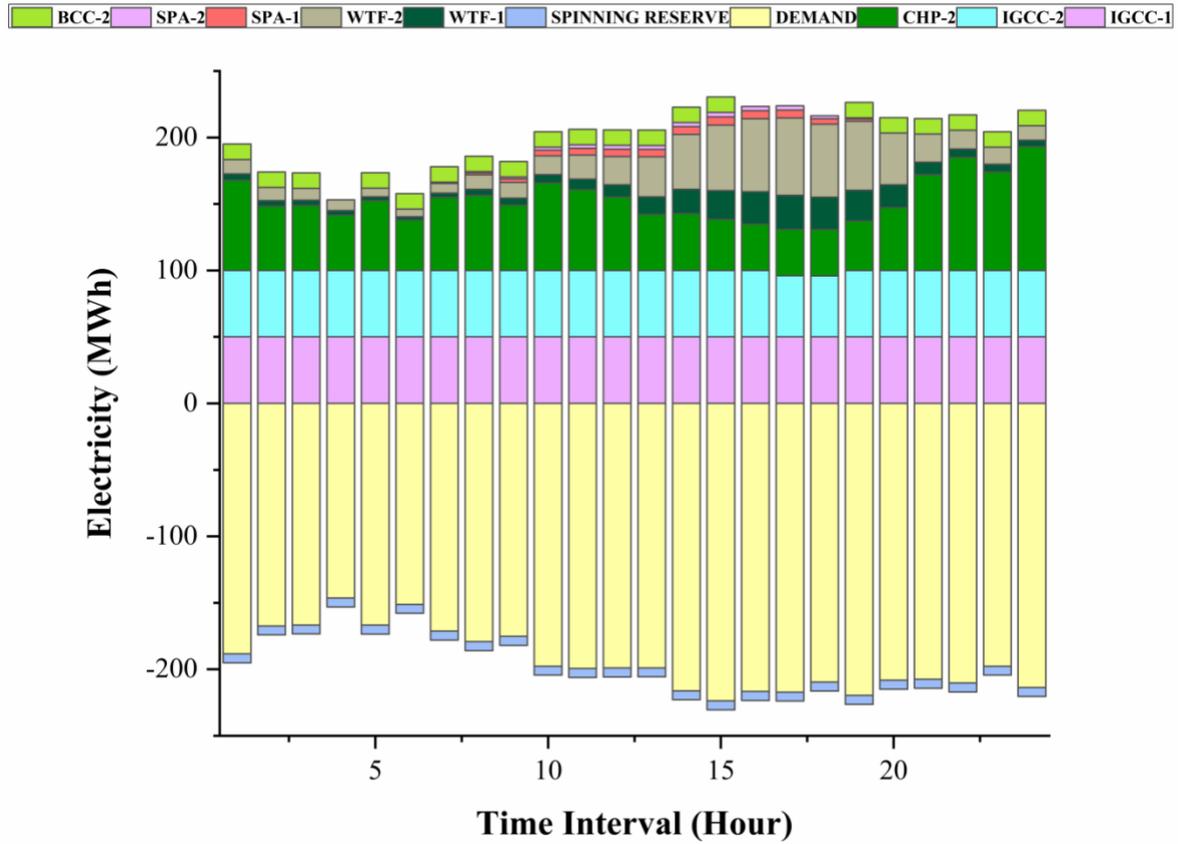

**Fig. 5.** Daily electricity plan for scenario *w31*.

Scenario *w31* includes summer data of hourly wind speed, solar radiation, and air temperature. As noticed in Fig. 6, the operating power of SPA-1 is much higher than SPA-2. It means that the panel efficiency of SPA-1 is higher than SPA-2, and mono-crystalline PV cells are more convenient for this weather profile than polycrystalline. Moreover, the operating powers of SPA-1 and SPA-2 are zero over the night due to zero solar radiation. They generate electricity at the intervals from 9 to 19. Both wind turbines produce much more electricity at time intervals between 14 and 20. Fig. 6 and Fig. 7 show hourly excess carbon dioxide emission and carbon dioxide limit over the day in scenario *w31,* respectively. As seen in Figs. 6 and 7, excess carbon dioxide is much lower than the carbon dioxide limit at the intervals from 15 to 19 due to the maximum electricity generation of WT-1 and WT-2. The highest excess carbon dioxide generation is seen at time interval 24 because all sustainable generators do not run at this interval and the electricity demand is the highest.

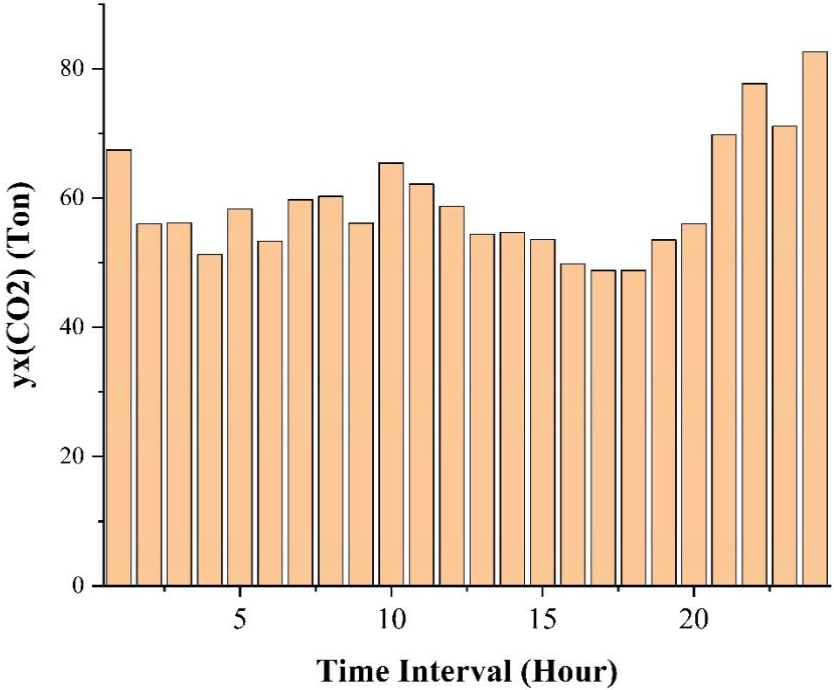

**Fig. 6.** Excess carbon dioxide emission over a day for scenario *w31*.

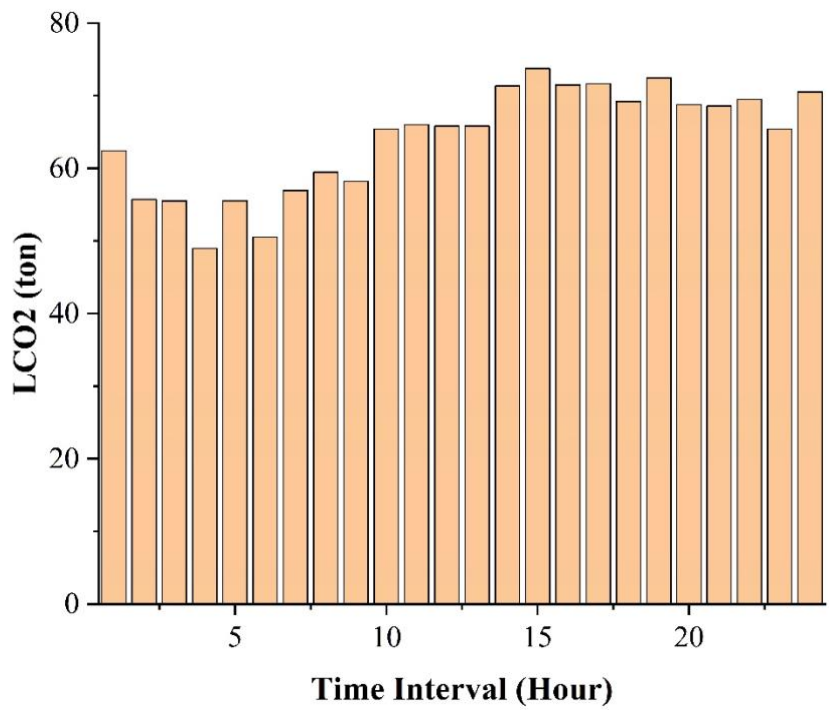

**Fig. 7.** Hourly carbon dioxide limit for scenario *w31*.

Fig. 8 shows hourly carbon trading cost or revenue over the day for scenario *w31*.

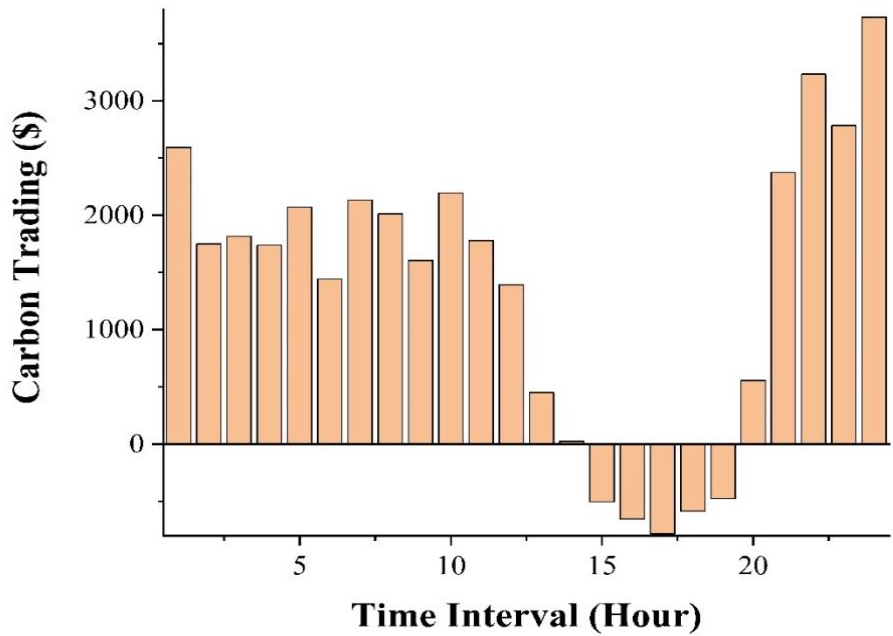

**Fig. 8.** Hourly carbon trading cost or revenue ($C^{CT}$) for scenario *w31*.

As seen in Fig. 8, the remaining carbon dioxide limit of the MEMG is sold at time intervals between 15 and 19. $C^{CT}$ is negative at these intervals. It means that excess carbon dioxide is lower than carbon dioxide limit. Excess carbon dioxide is higher than the carbon dioxide limit at the other intervals. Thus, extra carbon dioxide limit is bought by the MEMG at the other intervals and $C^{CT}$ is positive.

## 5. Conclusions

The aim of this study is to determine the optimal equipment selection for the design and scheduling of a Mixed Energy Management Grid (MEMG) under uncertainty, in order to meet the demands for electricity, heat, and natural gas. A novel two-stage Mixed-Integer Nonlinear Programming (MINLP) model is developed, considering thirty-two different scenarios. The uncertainties addressed in this study include carbon dioxide regulation policies, air temperature, wind speed, solar radiation, carbon dioxide trading prices or taxes, and natural gas prices. The candidate equipment under consideration includes wind turbine farms, solar panel arrays, biomass-fired generators, biomass combined cycles, Combined Heat and Power (CHP) units, conventional generators, electricity storage units, Integrated Gasification Combined Cycle (IGCC) systems, heat pumps, and Power-to-Gas (PTG) systems.

Key findings of the study are summarized as follows:

- The deterministic models favor wind turbine farms over biomass combined cycles and biomass-fired generators, whereas the stochastic models present different results.
- The number of scenarios considered in the two-stage stochastic models (whether sixteen with carbon trading, sixteen with carbon taxing, or the full thirty-two scenarios) does not significantly impact the optimal equipment selection.
- In the PTG-integrated MEMG, the model selects the maximum capacity of synthetic natural gas (SNG) generation, even though the natural gas demand is significantly lower than this amount.
- The two-stage optimization approach, when applied to the thirty-two scenarios, does not select wind turbines to be installed at their maximum rated capacity, even when they are explicitly required to be part of the grid.

**Acknowledgments**

This study has been performed benefiting from the 2232 International Fellowship for Outstanding Researchers Program of TUBITAK (Project No: 118C245). Still, the whole responsibility of this publication belongs to its authors. We also thank Prof. Metin Türkay for the valuable suggestions.